\centerline {\bf Four conjectures in Nonlinear Analysis}\par
\bigskip
\bigskip
\centerline {BIAGIO RICCERI}\par
\bigskip
\bigskip
{\bf Abstract.} In this chapter, I formulate four challenging conjectures in Nonlinear Analysis. More precisely: a conjecture on the
Monge-Amp\`ere equation; a conjecture on an eigenvalue problem; a conjecture on a non-local problem; a conjecture
on disconnectedness versus infinitely many solutions.\par
\bigskip
\bigskip
\bigskip
\bigskip
In this chapter, I intend to formulate four challenging conjectures in Nonlinear Analysis
which have their roots in certain results that I have obtained in the past years.\par
\bigskip
{\bf 1. A conjecture on the Monge-Amp\`ere equation}\par
\bigskip
CONJECTURE 1.1.  - {\it Let $\Omega\subset {\bf R}^n$ $(n\geq 2)$ be a non-empty open bounded set and let $h:\Omega\to {\bf R}$
be a non-negative continuous function.\par
Then, each $u\in C^2(\Omega)\cap C^1(\overline {\Omega})$ satisfying in $\Omega$ the Monge-Amp\`ere equation
$$\hbox {\rm det}(D^2u)=h$$
has the following property:
$$\nabla(\Omega)\subseteq \hbox {\rm conv}(\nabla(\partial\Omega))\ .$$}
\smallskip
This conjecture is motivated by [26] where I proved that it is true for $n=2$. I am going to produce such a proof here.\par
\smallskip
In what follows, $\Omega$ is a non-empty relatively compact and open set in a topological space $E$, with $\partial \Omega\neq \emptyset$, and
$Y$ is a real locally convex Hausdorff topological vector space. $\overline {\Omega}$ and $\partial \Omega$ denote the closure and the
boundary of $\Omega$, respectively. Since $\overline {\Omega}$ is compact, $\partial\Omega$, being closed, is compact too. 
\smallskip
Let us first recall some well-known definitions.\par
\smallskip
Let $S$ be a subset of $Y$ and let $y_0\in S$. As usual, we say that $S$ is supported at $y_0$ 
 if there exists $\varphi\in Y^*\setminus \{0\}$
such that $\varphi(y_0)\leq \varphi(y)$ for all $y\in S$. If this happens, of course
$y_0\in\partial S$.
\smallskip
Further, extending a maximum principle definition for real-valued functions,
a continuous function $f:\overline {\Omega}\to Y$ is said to  satisfy the convex hull property in $\overline {\Omega}$ (see [7], [13] and references therein)  if
$$f(\Omega)\subseteq \overline {\hbox {\rm conv}}(f(\partial \Omega))\ ,$$
$\overline {\hbox {\rm conv}}(f(\partial \Omega))$ being the closed convex hull of $f(\partial \Omega)$.
\smallskip
When dim$(Y)<\infty$, since $f(\partial\Omega)$ is compact, conv$(f(\partial\Omega))$ is compact too and so $\overline 
{\hbox {\rm conv}}(f(\partial \Omega))=\hbox {\rm conv}(f(\partial\Omega))$.\par
\smallskip
 A function $\psi:Y\to {\bf R}$ is said to be quasi-convex if, for each $r\in {\bf R}$, the set $\psi^{-1}(]-\infty,r])$ is convex.\par
\smallskip
Notice the following proposition:\par
\medskip
PROPOSITION 1.1. - {\it For each pair $A, B$ of non-empty subsets of $Y$, the following assertions are equivalent:\par
\noindent
$(a_1)$\hskip 5pt $A\subseteq \overline{\hbox {\rm conv}}(B)\ .$\par
\noindent
$(a_2)$\hskip 5pt For every continuous and quasi-convex function $\psi:Y\to {\bf R}$, one has
$$\sup_A\psi\leq \sup_B\psi\ .$$}
\smallskip
PROOF. Let $(a_1)$ hold. Fix any continuous and quasi-convex function $\psi:Y\to {\bf R}$. Fix $\tilde y\in A$.
Then, there is a net $\{y_{\alpha}\}$ in conv$(B)$ converging to $\tilde y$. So, for each $\alpha$, we have
$y_{\alpha}=\sum_{i=1}^k\lambda_i z_i$, where $z_i\in B$, $\lambda_i\in [0,1]$ and $\sum_{i=1}^k\lambda_i=1$.
By quasi-convexity, we have
$$\psi(y_{\alpha})=\psi\left ( \sum_{i=1}^k\lambda_iz_i\right )\leq \max_{1\leq i\leq k}\psi(z_i)\leq \sup_B\psi$$
and so, by continuity,
$$\psi(\tilde y)=\lim_{\alpha}\psi(y_{\alpha})\leq\sup_B\psi$$
which yields $(a_2)$.\par
Now, let $(a_2)$ hold. Let $x_0\in A$. If $x_0\not\in \overline{\hbox {\rm conv}}(B)$, by the standard separation theorem, there
would be $\psi\in Y^*\setminus\{0\}$ such that $\sup_{\overline{\hbox {\rm conv}}(B)}\psi<\psi(x_0)$, against $(a_2)$. So,
$(a_1)$ holds.\hfill $\bigtriangleup$\par
\medskip
Clearly, applying Proposition 1.1, we obtain the following one:\par
\medskip
PROPOSITION 1.2. - {\it For any continuous function $f:\overline {\Omega}\to Y$, the following assertions are equivalent:\par
\noindent
$(b_1)$\hskip 5pt $f$ satisfies the convex hull property in $\overline {\Omega}$\ .\par
\noindent
$(b_2)$\hskip 5pt For every continuous and quasi-convex function $\psi:Y\to {\bf R}$, one has
$$\sup_{x\in \Omega}\psi(f(x))=\sup_{x\in\partial \Omega}\psi(f(x))\ .$$ }\par
\medskip
In view of Proposition 1.2, we now introduce the notion of convex hull-like property for functions defined in $\Omega$ only.\par
\medskip
DEFINITION 1.1. - A continuous function $f:\Omega\to Y$ is said to satisfy the convex hull-like property in $\Omega$ if, for every continuous and quasi-convex function
$\psi:Y\to {\bf R}$, there exists $x^*\in\partial \Omega$ such that
$$\limsup_{x\to x^*}\psi(f(x))=\sup_{x\in \Omega}\psi(f(x))\ .$$
\medskip
We have\par
\medskip
PROPOSITION 1.3. - {\it Let $g:\overline {\Omega}\to Y$ be a continuous function and let $f=g_{|\Omega}$.\par
Then, the following assertions are equivalent:\par
\noindent
$(c_1)$\hskip 5pt $f$ satisfies the convex hull-like property in $\Omega$\ .\par
\noindent
$(c_2)$\hskip 5pt $g$ satisfies the convex hull property in $\overline {\Omega}$\ .}\par
\medskip
PROOF. Let $(c_1)$ hold. Let $\psi:Y\to {\bf R}$ be any continuous and quasi-convex function.
Then, by Definition 1.1, there exists $x^*\in\partial \Omega$ such that
$$\limsup_{x\to x^*}\psi(f(x))=\sup_{x\in \Omega}\psi(f(x))\ .$$
But
$$\limsup_{x\to x^*}\psi(f(x))=\psi(g(x^*))$$
and hence
$$\sup_{x\in\partial \Omega}\psi(g(x))=\sup_{x\in \Omega}\psi(g(x))\ .$$
So, by Proposition 1.2, $(c_2)$ holds.\par
Now, let $(c_2)$ hold. Let $\psi:Y\to {\bf R}$ be any continuous and quasi-convex function. Then, by
Proposition 1.2, one has
$$\sup_{x\in\partial \Omega}\psi(g(x))=\sup_{x\in \Omega}\psi(g(x))\ .$$
Since $\partial \Omega$ is compact and $\psi\circ g$ is continuous,
there exists $x^*\in\partial \Omega$ such that 
$$\psi(g(x^*))=\sup_{x\in\partial \Omega}\psi(g(x))\ .$$
But
$$\psi(g(x^*))=\lim_{x\to x^*}\psi(f(x))$$
and, by continuity again,
$$\sup_{x\in \Omega}\psi(g(x))=\sup_{x\in \overline {\Omega}}\psi(g(x))$$
and so 
$$\lim_{x\to x^*}\psi(f(x))=\sup_{x\in \Omega}\psi(f(x))$$
which yields $(c_1)$.\hfill $\bigtriangleup$\par
\medskip
The central result is as follows:\par
\medskip
THEOREM 1.1. - {\it For any continuous function
 $f:\Omega\to Y$, at least one of the following assertions holds:\par
\noindent
$(i)$\hskip 5pt $f$ satisfies the convex hull-like property in $\Omega$\ .\par
\noindent
$(ii)$\hskip 5pt There exists a non-empty open set $X\subseteq \Omega$, with $\overline {X}\subseteq \Omega$,  satisfying the following property:
for every continuous function $g:\Omega\to Y$, there exists $\tilde\lambda\geq 0$
such that, for each $\lambda>\tilde\lambda$, the set $(g+\lambda f)(X)$
is supported at one of its points.}\par
\smallskip
PROOF.  Assume that $(i)$ does not hold. So, we are assuming that there exists
a continuous and quasi-convex function $\psi:Y\to {\bf R}$ such that
$$\limsup_{x\to z}\psi(f(x))<\sup_{x\in \Omega}\psi(f(x)) \eqno{(1.1)}$$
for all $z\in \partial \Omega$.\par
In view of $(1.1)$, for each $z\in \partial \Omega$, there exists 
an open neighbourhood $U_z$ of $z$ such that
$$\sup_{x\in U_z\cap \Omega}\psi(f(x))<\sup_{x\in \Omega}\psi(f(x))\ .$$
Since $\partial\Omega$ is compact, there are finitely many $z_1,...,z_k\in
\partial \Omega$ such that
$$\partial \Omega\subseteq \bigcup_{i=1}^kU_{z_i}\ .\eqno{(1.2)}$$
Put
$$U=\bigcup_{i=1}^kU_{z_i}\ .$$
Hence
$$\sup_{x\in U\cap \Omega}\psi(f(x))=\max_{1\leq i\leq k}\sup_{x\in U_{z_i}\cap \Omega}\psi(f(x))<
\sup_{x\in \Omega}\psi(f(x))\ .$$
Now, fix a number $r$ so that
$$\sup_{x\in U\cap \Omega}\psi(f(x))<r<\sup_{x\in \Omega}\psi(f(x)) \eqno{(1.3)}$$
and set
$$K=\{x\in \Omega : \psi(f(x))\geq r\}\ .$$
Since $f, \psi$ are continuous, $K$ is closed in $\Omega$. But,
since $K\cap U=\emptyset$ and $U$ is open, in view of $(1.2)$, 
$K$ is closed in $E$. Hence, $K$ is compact since  $\overline {\Omega}$ is so.
By $(1.3)$, we can fix $\bar x\in \Omega$ such that $\psi(f(\bar x))>r$. Notice that the set
$\psi^{-1}(]-\infty,r])$ is closed and convex. So, thanks to
the standard separation theorem,  there exists
a non-zero continuous linear functional $\varphi:Y\to {\bf R}$ such that
$$\varphi(f(\bar x))<\inf_{y\in \psi^{-1}(]-\infty,r])}\varphi(y)\ .\eqno{(1.4)}$$
Then, from $(1.4)$, it follows
$$\varphi(f(\bar x))<\inf_{x\in \Omega\setminus K}\varphi(f(x))\ .$$
Now, choose $\rho$ so that
$$\varphi(f(\bar x))<\rho<\inf_{x\in \Omega\setminus K}\varphi(f(x))$$
and set
$$X=\{x\in \Omega : \varphi(f(x))<\rho\}\ .$$
Clearly, $X$ is a non-empty open set contained in $K$. Now, let $g:\Omega\to Y$ be any continuous function. 
Set
$$\tilde\lambda=\inf_{x\in X}{{\varphi(g(x))-\inf_{z\in K}\varphi(g(z))}\over {\rho-\varphi(f(x))}}\ .$$
Fix $\lambda>\tilde\lambda$. So, there is $x_0\in X$ such that
$${{\varphi(g(x_0))-\inf_{z\in K}\varphi(g(z))}\over {\rho-\varphi(f(x_0))}}<\lambda\ .$$
  From this, we get
  $$\varphi(g(x_0))+\lambda\varphi(f(x_0))<\lambda \rho+\inf_{z\in K}\varphi(g(z))\ .\eqno{(1.5)}$$
By continuity and compactness, there 
exists $\hat x\in K$ such that
$$\varphi(g(\hat x)+\lambda f(\hat x))\leq\varphi(g(x))+\lambda f(x))
\eqno{(1.6)}$$
for all $x\in K$.
Let us prove that $\hat x\in X$.
Arguing by contradiction, assume that $\varphi(f(\hat x))\geq \rho$.
Then, taking $(1.5)$ into account, we would have
$$\varphi(g(x_0))+\lambda\varphi(f(x_0))<\lambda\varphi(f(\hat x))+\varphi(g(\hat x))$$
contradicting $(1.6)$. So, it is true that $\hat x\in X$, and, by $(1.6)$, the set
$(g+\lambda f)(X)$ is supported at its point $g(\hat x)+\lambda f(\hat x)$.
\hfill $\bigtriangleup$\par
\medskip
An application of Theorem 1.1 shows a strongly bifurcating behaviour of certain equations in ${\bf R}^n$.\par
\medskip
THEOREM 1.2. - {\it Let $\Omega$ be a non-empty bounded open subset of ${\bf R}^n$ and let $f:\Omega\to {\bf R}^n$ a
continuous function.\par
Then,  at least one of the following assertions holds:\par
\noindent
$(d_1)$\hskip 5pt $f$ satisfies the convex hull-like property in $\Omega$\ .\par
\noindent
$(d_2)$\hskip 5pt There exists a non-empty open set $X\subseteq \Omega$, with $\overline {X}\subseteq \Omega$, satisfying the following property:
for every continuous function $g:\Omega\to {\bf R}^n$, there exists $\tilde\lambda\geq 0$
such that, for each $\lambda>\tilde\lambda$, there exist $\hat x\in X$ and two
sequences $\{y_k\}$, $\{z_k\}$ in ${\bf R}^n$, with
$$\lim_{k\to \infty}y_k=\lim_{k\to \infty}z_k=g(\hat x)+\lambda f(\hat x)\ ,$$
such that, for each $k\in {\bf N}$, one has\par
\noindent
$(j)$\hskip 5pt the equation
$$g(x)+\lambda f(x)=y_k$$
has no solution in $X$\ ;\par
\noindent
$(jj)$\hskip 5pt the equation
$$g(x)+\lambda f(x)=z_k$$
has two distinct solutions $u_k, v_k$ in $X$ such that
$$\lim_{k\to \infty}u_k=\lim_{k\to \infty}v_k=\hat x\ .$$}\par
\smallskip
PROOF.  Apply Theorem 1.1 with $E=Y={\bf R}^n$. Assume that $(d_1)$ does not hold.  Let $X\subseteq \Omega$ be an open
set as in $(ii)$ of Theorem 1.1. Fix any continuous function $g:\Omega\to {\bf R}^n$. 
Then, there is
some $\tilde\lambda\geq 0$ such that, for each $\lambda>\tilde\lambda$, there exists $\hat x\in X$
such that the set $(g+\lambda f)(X)$ is supported at $g(\hat x)+\lambda f(\hat x)$. As we observed at the
beginning, this implies that
$g(\hat x)+\lambda f(\hat x)$ lies in the boundary of $(g+\lambda f)(X)$.
Therefore, we can find a sequence $\{y_k\}$ in ${\bf R}^n\setminus (g+\lambda f)(X)$ converging
to $g(\hat x)+\lambda f(\hat x)$. So, such a sequence satisfies $(j)$. For each $k\in {\bf N}$,
denote by $B_k$ the open ball of radius ${{1}\over {k}}$ centered at $\hat x$. Let $k$ be such that
$B_k\subseteq X$. The set $(g+\lambda f)(B_k)$ is not open since  its boundary contains the point
$g(\hat x)+\lambda f(\hat x)$. Consequently, by the invariance of domain theorem ([29], p. 705), 
the function $g+\lambda f$ is not injective in $B_k$. So, there are $u_k, v_k\in B_k$, with $u_k\neq v_k$
such that
$$g(u_k)+\lambda f(u_k)=g(v_k)+\lambda f(v_k)\ .$$
Hence, if we take
$$z_k=g(u_k)+\lambda f(u_k)\ ,$$
the sequences $\{u_k\}, \{v_k\}, \{z_k\}$ satisfy $(jj)$ and the proof is complete.\hfill $\bigtriangleup$
\par
\medskip
REMARK 1.1. - Notice that, in general, Theorem 1.2 is no longer true when $f:\Omega\to {\bf R}^m$ with $m>n$. In this connection,
consider the case $n=1$, $m=2$, $\Omega=]0,\pi[$ and $f(\theta)=(\cos\theta,\sin\theta)$ for $\theta\in [0,\pi]$. So,
for each $\lambda>0$, on the one hand, the function $\lambda f$ is injective, while, on the other hand, $\lambda f(]0,\pi[)$
is not contained in conv$(\{f(0),f(\pi)\})$.\par
\medskip
If $S\subseteq {\bf R}^n$ is a non-empty open set, $x\in S$ and $h:S\to {\bf R}^n$ is a $C^1$ function, we denote
by det$(J_h(x))$ the Jacobian determinant of $h$ at $x$.\par
\smallskip
A very recent and important result by J. Saint Raymond ([27]) states what follows (for anything concerning the topological dimension we refer to [8]):\par
\medskip
THEOREM 1.A ([27], Theorem 10). - {\it  Let $A\subseteq {\bf R}^n$ be a non-empty open set and $\varphi:A\to {\bf R}^n$  a $C^1$  function such that
the topological dimension of the set 
$$\{x\in A : \hbox {\rm det}(J_{\varphi}(x))=0\}$$
is not positive.\par
Then, the function $\varphi$ is open.}\par
\medskip
A joint application of Theorem 1.1 and Theorem 1.A gives
\medskip
THEOREM 1.3. - {\it  Let
 $f:\Omega\to {\bf R}^n$ be a $C^1$ function.\par
Then,  at least one of the following assertions holds:\par
\noindent
$(a_1)$\hskip 5pt $f$ satisfies the convex hull-like property in $\Omega$\ .\par
\noindent
$(a_2)$\hskip 5pt There exists a non-empty open set $X\subseteq \Omega$, with $\overline {X}\subseteq \Omega$, satisfying the following property:
for every continuous  function $g:\Omega\to {\bf R}^n$ which is $C^1$ in $X$, there exists $\tilde\lambda\geq 0$ such that, for each 
$\lambda>\tilde\lambda$, the topological dimension of the 
set 
$$\{x\in X : \hbox {\rm det}(J_{g+\lambda f}(x))=0\}$$
 is greater than or equal $1$.}\par
\smallskip
PROOF.  Assume that $(a_1)$ does not hold. Let $X$ be an open set as in $(ii)$ of Theorem 1.1. Let $g:\Omega\to {\bf R}^n$ be a continuous function which
is $C^1$ in $X$. Then, there is
some $\tilde\lambda\geq 0$ such that, for each $\lambda>\tilde\lambda$, there exists $\hat x\in X$
such that the set $(g+\lambda f)(X)$ is supported at $g(\hat x)+\lambda f(\hat x)$.  As already remarked, this implies that
$g(\hat x)+\lambda f(\hat x)\in \partial(g+\lambda f)(X)$ and so $(g+\lambda f)(X)$ is not open. Now, $(a_2)$ is a direct consequence
of Theorem 1.A.\hfill $\bigtriangleup$\par
\medskip
In turn, here is a consequence of Theorem 1.3 when $n=2$.\par
\medskip
THEOREM 1.4. - {\it Let $\Omega$ be a non-empty bounded open set of ${\bf R}^2$, let
$h:\Omega\to {\bf R}$ be a continuous function and let $\alpha, \beta:\Omega\to {\bf R}$ be two $C^1$ functions such that
$|\alpha_x\beta_y-\alpha_y\beta_x|+|h|>0$ and $(\alpha_x\beta_y-\alpha_y\beta_x)h\geq 0$ in $\Omega$.\par
Then, any $C ^1$ solution $(u,v)$ in $\Omega$ of the system
$$\cases{u_xv_y-u_yv_x=h\cr & \cr 
\beta_y u_x-\beta_x u_y-\alpha_yv_x+\alpha_xv_y =0\cr}\eqno{(1.7)}$$
 satisfies the convex hull-like property in $\Omega$.}\par
\smallskip
PROOF. Arguing by contradiction, assume that $(u,v)$ does not satisfy the convex hull-like property in $\Omega$. Then, by Theorem 1.3, applied taking
$f=(u,v)$ and $g=(\alpha,\beta)$, there exist $\lambda>0$ and $(\hat x,\hat y)\in \Omega$ such that
$$\hbox {\rm det}(J_{g+\lambda f}(\hat x,\hat y))=0\ .$$
On the other hand, for each $(x,y)\in \Omega$, we have
$$\hbox {\rm det}(J_{g+\lambda f}(x, y))=(u_xv_y-u_yv_x)(x,y)\lambda^2+(\beta_y u_x-\beta_x u_y-\alpha_yv_x+\alpha_xv_y)(x,y)\lambda+
(\alpha_x\beta_y-\alpha_y\beta_x)(x,y)$$
and hence
$$h(\hat x,\hat y)\lambda^2+(\alpha_x\beta_y-\alpha_y\beta_x)(\hat x,\hat y)=0$$
which is impossible in view of our assumptions.\hfill $\bigtriangleup$\par
\medskip
Finally, taking Proposition 1.3 in mind, here is the proof of Conjecture 1.1 when $n=2$: \par
\medskip
THEOREM 1.5. - {\it Let $\Omega$ be a non-empty bounded open subset of ${\bf R}^2$, let $h:\Omega\to {\bf R}$ be a continuous non-negative function and
let $w\in C^2(\Omega)$ be a function satisfying in $\Omega$ the  Monge-Amp\`ere equation
$$w_{xx}w_{yy}-w_{xy}^2=h\ .$$
Then, the gradient of $w$ satisfies the convex hull-like property in $\Omega$.}\par
\smallskip
PROOF. It is enough to observe that $(w_x,w_y)$ is a $C^1$ solution in $\Omega$ of the system $(1.7)$ with
$\alpha(x,y)=-y$ and $\beta(x,y)=x$ and that such $\alpha, \beta$ satisfy the assumptions of Theorem 1.4.\hfill
$\bigtriangleup$
\par
\bigskip
{\bf 2. A conjecture on an eigenvalue problem}\par
\bigskip
CONJECTURE 2.1. - {\it Let $n\geq 2$ and let $\Omega =\{x\in {\bf R}^n
:a<|x|<b\}$, with $0<a<b$.\par
Then, there exists $\lambda>0$ such that the problem
$$\cases {\Delta u=\lambda\sin u & in
$\Omega$\cr & \cr u=0 & on
$\partial \Omega$\cr}$$
has at least one non-zero classical solution.}\par
\medskip
The above conjecture has its roots in Pohozaev identity ([19]). Let me recall it.\par
\smallskip
So, let $\Omega\subset {\bf R}^n$ be a smooth bounded domain, and let $f:{\bf R}\to {\bf R}$ be
a continuous function. Put
$$F(\xi)=\int_0^{\xi} f(t)dt$$
for all $\xi\in {\bf R}$. For $\lambda>0$, consider the problem
 $$\cases {-\Delta u=\lambda f(u) & in $\Omega$\cr & \cr u=0 & on $\partial
\Omega$\ .\cr} \eqno{(P_{\lambda f})}$$
In the sequel, a classical solution of problem $(P_{\lambda f})$ is any
$u\in C^2(\Omega)\cap C^1(\overline {\Omega})$, zero on $\partial \Omega$, satisfying
the equation pointwise in $\Omega$. Set
$$\Lambda_f = \{\lambda>0 : (P_{\lambda f})\hskip 5pt \hbox {\rm has a non-zero classical solution}\}\ .$$
When $n\geq 2$, the Pohozaev identity tells us
that, if $u$ is a classical solution of $(P_{\lambda f})$, then
one has
$${{2-n}\over {2}}\int_{\Omega}|\nabla u(x)|^2dx
+n\lambda\int_{\Omega}F(u(x))dx=
{{1}\over {2}}\int_{\partial \Omega}|\nabla u(x)|^2 x\cdot
\nu(x) ds \eqno{(2.1)}$$
where $\nu$ denotes the unit outward normal to $\partial \Omega$.
\par
\smallskip
  From $(2.1)$, in particular, it follows that,
 if $\Omega$ is star-shaped with respect to $0$
 (so $x\cdot \nu(x)\geq 0$
  on $\partial \Omega$), then the set
  $\Lambda_f$ is empty in the two following cases:\par
\smallskip
\noindent
$(a)$\hskip 5pt $f(\xi)=|\xi|^{p-2}\xi$ with $n\geq 3$ and $p\geq {{2n}\over {n-2}}$\ ;\par
\smallskip
\noindent
$(b)$\hskip 5pt $\sup_{\xi\in {\bf R}}F(\xi)=0$\ .\par
\smallskip
A natural question arises: what about problem $(P_{\lambda f})$ in cases $(a)$ and $(b)$ when $\Omega$ is not star-shaped ?\par
\smallskip
It is very surprising to realize that, while a great amount of research has been produced on case $(a)$ (see, for instance, [1]-[5], [12], [14], [17], [18]), 
apparently the only papers dealing with case $(b)$ are [9]-[11], [23]. \par
\smallskip
In [11], the following result has been pointed out:\par
\medskip
THEOREM 2.1. - {\it  
Let $n\geq 2$ and $\Omega =\{x\in {\bf R}^n
:a<|x|<b\}$ with $0<a<b$.\par
 Then for every continuous function $f:{\bf R}\to {\bf R}$, 
with $\sup_{\xi\in {\bf R}}F(\xi)=0$, and every $\lambda >0$, problem $(P_{\lambda f})$ has no
radially symmetric non-zero classical solutions.}\par
\smallskip
PROOF. 
Let $f:{\bf R}\to {\bf R}$ be a continuous function,
with $\sup_{\xi\in {\bf R}}F(\xi)=0$, let $\lambda >0$, and let $u$
 be a radially symmetric classical solution of $(P_{\lambda f})$. Then
$$\cases {
-(r^{n-1}u^{\prime }(r))^{\prime}=\lambda r^{n-1}f(u(r)) & for $r\in
(a,b)$ \cr & \cr
u(a)=u(b)=0\ ,\cr}$$
that is
$$\cases {
u^{\prime \prime }(r)+{{n-1}\over {r}}u^{\prime }(r)+\lambda f(u(r))=0 &
for $r\in (a,b)$ \cr & \cr
u(a)=u(b)=0\ .\cr}\eqno{(2.2)}$$
Multiplying both sides of the equation in (2.2) by $u^{\prime }$, we have
$$
u^{\prime \prime }(r)u^{\prime }(r)+{{n-1}\over {r}}(u^{\prime
}(r))^{2}+\lambda f(u(r))u^{\prime }(r)=0  \eqno{(2.3)}$$
for all $r\in (a,b)$.
Let $r_{1}\in (a,b)$ be such that $u^{\prime }(r_{1})=0$. Define
$$I_{r_{1}}(r)={{1}\over {2}}\left\vert u^{\prime }(r)\right\vert
^{2}+(n-1)\int_{r_{1}}^{r}{{(u^{\prime }(t))^{2}}\over {t}}dt+\lambda F(u(r))$$
for all $r\in [a,b]$.
Then (2.3) shows that $I_{r_{1}}^{\prime }(r)=0$ for all $r\in (a,b)$ and so, for some $c\in {\bf R}$, one has
$$
I_{r_{1}}(r)=c$$
for all $r\in [a,b]$. Since
$$
I_{r_{1}}(r_{1})=0+0+\lambda F(u(r_{1}))\leq 0\ ,$$
we have $c\leq 0$. On the other hand, since
$$
I_{r_{1}}(b)={{1}\over {2}}|u^{\prime }(b)|
^{2}+(n-1)\int_{r_{1}}^{b}{{(u^{\prime }(t))^{2}}\over {t}}dt+0\geq 0\ ,
$$
have $c\geq 0$, and so $c=0$. In particular $I_{r_{1}}(b)=0$, which implies
$u^{\prime }(b)=0$, and consequently $u(r)=0$ for all $r\in [a,b]$,
as claimed.\hfill $\bigtriangleup$\par
\medskip
REMARK 2.1. - It is important to note the drastic difference 
 between cases $(a)$ and $(b)$ enlighted by Theorem 2.1 when $\Omega$ is an annulus. Actually, in this case, it was remarked in [13] that the problem
$$\cases {-\Delta u=\lambda |u|^{p-2}u & in $\Omega$\cr & \cr u=0 & on $\partial
\Omega$\cr}$$
has radially symmetric non-zero classical solutions for $p\geq {{2n}\over {n-2}}$ 
($n\geq 3$), and $\lambda>0$.\par
\medskip
Now, I recall a very general result proved in [23].\par
\smallskip
For any real Hilbert space
$X$, denote by ${\cal A}_{X}$ the set of all $C^1$
functionals $I:X\to {\bf R}$ such that $0$ is a global maximum of $I$
and $I'$ is Lipschitzian with Lipschitz constant less than $1$.
 Set
$$\gamma_X=\inf_{I\in {\cal A}_X}\inf\{\lambda>0 : x=\lambda I'(x)
\hskip 5pt \hbox {\rm for some} \hskip 5pt x\neq 0\}\ .$$
We have:\par
\medskip
THEOREM 2.2. - {\it For any real Hilbert space $(X,\langle\cdot,\cdot\rangle)$, 
with $X\neq \{0\}$, one has
$$\gamma_X=3\ .$$}\par
\medskip
We first prove 
\medskip
PROPOSITION 2.1. - {\it One has
$$\gamma_{\bf R}=3\ .$$}
\smallskip
PROOF. Let $I_0\in {\cal A}_{\bf R}$ and let $L<1$ be the Lipschitz constant
of $I_0'$. Set
$$I=I_0-I_0(0)\ .$$
Fix $\lambda\in ]0,3]$. Let us prove that $0$ is
the only solution of the equation
$$x=\lambda I'(x)\ .$$
Arguing by contradiction, assume that
$$x_0=\lambda I'(x_0)$$
for some $x_0\neq 0$. It is not restrictive to assume that $x_0>0$
(otherwise, we would work with $I'(-x)$). Consider now the function
$g:{\bf R}\to {\bf R}$ defined by
$$g(x)=\cases {-{{x^2}\over {2}} & if $x<{{x_0}\over {3}}$ \cr & \cr
{{x^2}\over {2}}-{{2x_0x}\over {3}}+{{x_0^2}\over {9}} & if
${{x_0}\over {3}}\leq x\leq x_0$ \cr & \cr -{{x^2}\over {2}}+
{{4x_0x}\over {3}}-{{8x_0^2}\over {9}} & if $x_0>x$\ . \cr}$$
Clearly, $g\in C^1({\bf R})$. Let $x>0$ with $x\neq x_0$. Let us prove
that 
$$g'(x)<I'(x)\ .$$ We distinguish two cases. If $0<x\leq {{x_0}\over {3}}$,
We have
$$g'(x)=-x<-Lx\leq I'(x)\ .$$
If $x>{{x_0}\over {3}}$, We have
$$g'(x)={{x_0}\over {3}}-|x-x_0|<{{x_0}\over {3}}-L|x-x_0|={{\lambda I'(x_0)}
\over {3}}-L|x-x_0|\leq I'(x_0)-L|x-x_0|\leq
I'(x)\ . $$
So, in particular, we get
$$I\left ( {{4x_0}\over {3}}\right )=\int_0^{4x_0\over 3}I'(x)dx>\int_0^{4x_0\over 3}g'(x)dx=
g\left ( {{4x_0}\over {3}}\right )=0$$
which contradicts the fact that the function $I$ is non-positive, since $0$ is a global
maximum of $I_0$. From what we have just proven, it clearly follows that
$$3\leq \gamma_{\bf R}\ .$$ Now, fix any $\mu>1$. Continue to consider the
function $g$ defined above (for a fixed $x_0>0$). Clearly, the function
${{1}\over {\mu}}g$ belongs to ${\cal A}_{\bf R}$ and
$$x_0=3\mu{{g(x_0)}\over {\mu}}\ .$$
Of course, from this we infer that
$$\gamma_{\bf R}\leq 3\mu$$
and the conclusion clearly follows.\hfill $\bigtriangleup$\par
\medskip
{\it Proof of Theorem 2.2.}
First, let us prove that
$$\gamma_X\leq 3\ .\eqno{(2.4)}$$
To this end,
fix any $\varphi\in {\cal A}_{\bf R}$ and any $\lambda>0$ such that
$$t=\lambda \varphi'(t)$$
for some $t\neq 0$. Fix also $u\in X$, with $\|u\|=1$, and consider
the functional $I$ defined by putting
$$I(x)=\varphi(\langle u,x\rangle)$$
for all $x\in X$. It is readily seen that $I\in {\cal A}_X$. In particular,
note that
$$I'(x)=\varphi'(\langle u,x\rangle)u\ .$$
Finally, set
$$\hat x=\lambda\varphi(t) u\ .$$
Of course, $\hat x\neq 0$. Since
$$\langle u,\hat x\rangle =\lambda\varphi'(t)$$
we also have
$$\langle u,\hat x\rangle =t$$
and so
$$\hat x=\lambda I'(\hat x)\ .$$
 From this, it clearly follows that
$$\gamma_X\leq \gamma_{\bf R}$$
and so $(2.4)$ follows now from Proposition 2.1.\par
Now, let us prove that
$$3\leq\gamma_X\ .\eqno{(2.5)}$$
To this end, fix $I\in {\cal A}_X$, $\lambda>0$ and $x\in X\setminus \{0\}$ such that
$$x=\lambda I'(x)\ .\eqno{(2.6)}$$
Then, consider the function $\varphi:{\bf R}\to {\bf R}$ defined by
$$\varphi(t)=I\left ( {{tx}\over {\|x\|}}\right )$$
for all $t\in {\bf R}$. Clearly, $0$ is a global maximum for $\varphi$. Moreover,
$\varphi\in C^1({\bf R})$ and one has
$$\varphi'(t)=\left \langle I'\left ( {{tx}\over {\|x\|}}\right ), {{x}\over {\|x\|}}\right
\rangle\ .$$
Therefore, if $L$ is the Lipschitz constant of $I'$,
for each $t,s\in {\bf R}$, we have
$$|\varphi'(t)-\varphi'(s)|=\left | \left \langle 
I'\left ( {{tx}\over {\|x\|}}\right )-I'\left ( {{sx}\over {\|x\|}}\right ), 
{{x}\over {\|x\|}}\right \rangle\right |$$
$$\leq \left \| I'\left ( {{tx}\over {\|x\|}}\right )-I'\left ( {{sx}\over {\|x\|}}\right )\right \|\leq
L|t-s|\ .$$
This shows that $\varphi'$ is a contraction,
and so $\varphi\in {\cal A}_{\bf R}$. Now, from $(2.6)$,
we get
$$\|x\|=\lambda \left \langle I'(x),{{x}\over {\|x\|}}\right \rangle$$
that is
$$\|x\|=\lambda \varphi'(\|x\|)\ .$$
 From this, we infer that 
$$\gamma_{\bf R}\leq \gamma_X\ .$$
So $(2.5)$ follows from Proposition 2.1, and the proof is complete.\hfill $\bigtriangleup$\par
\medskip
Now, for each $L>0$, denote by ${\cal C}_L$ the class of all Lipschitzian functions $f:{\bf R}\to {\bf R}$, with Lipschitz constant $L$, such that $f(0)=0$ and $\sup_{\xi\in {\bf R}}F(\xi)=0$. Also
denote by
$\lambda_1$ the first eigenvalue of the problem
$$\cases {-\Delta u=\lambda u & in $\Omega$\cr & \cr
u=0 & on $\partial \Omega$\ .\cr}$$
\smallskip
From Theorem 2.2, it directly follows that
$$\inf_{f\in {\cal C}_L}\inf\Lambda_f\geq {{3\lambda_1}\over {L}}\ .$$
\smallskip
In [10],  X. L. Fan obtained the finer inequality
$$\inf_{f\in {\cal C}_L}\inf\Lambda_f>{{3\lambda_1}\over {L}}\ .$$
\smallskip
Conjecture 2.1 says that $\Lambda_f\neq\emptyset$ for $f(\xi)=-\sin\xi$, $\Omega$ being an annulus. Due to what precedes,
 if Conjecture 2.1 is true, then $\lambda$ must necessarily be larger than $3\lambda_1$.\par
\vfill\eject
{\bf 3. A conjecture on a non-local problem}\par
\bigskip
CONJECTURE 3.1. - {\it Let $a\geq 0$, $b>0$ and let $\Omega\subset {\bf R}^n$ be a smooth bounded domain, with $n> 4$.\par
Then, for each $\lambda>0$ large enough and for each convex set $C\subseteq L^2(\Omega)$ whose closure in $L^2(\Omega)$
contains $H^1_0(\Omega)$, there exists $v^*\in C$ such that the problem
$$\cases {-\left ( a+b\int_{\Omega}|\nabla u(x)|^2dx\right )\Delta u =|u|^{{4}\over {n-2}}u+\lambda(u-v^*(x)) & in $\Omega$\cr
& \cr u=0 & on $\partial\Omega$\cr}$$
has at least three weak solutions, two of which are global minima in $H^1_0(\Omega)$ of the functional
$$u\to {{a}\over {2}}\int_{\Omega}|\nabla u(x)|^2dx+{{b}\over {4}}\left ( \int_{\Omega}|\nabla u(x)|^2dx\right ) ^2-{{n-2}\over
{2n}}\int_{\Omega}|u(x)|^{{2n}\over {n-2}}dx
-{{\lambda}\over {2}}\int_{\Omega}|u(x)-v^*(x)|^2dx\ .$$}\par
\medskip
Conjecture 3.1 comes from the results I have obtained in [25]. I am going to reproduce them here.\par
\smallskip
Let $a, b, \Omega$ be as in Conjecture 3.1.\par
\smallskip
On the Sobolev space $H^1_0(\Omega)$, we consider the scalar product
$$\langle u,v\rangle=\int_{\Omega}\nabla u(x)\nabla v(x)dx$$
and the induced norm
$$\|u\|=\left ( \int_{\Omega}|\nabla u(x)|^2dx\right )^{1\over 2}\ .$$
Denote by ${\cal A}$ the class of all Carath\'eodory functions $f:\Omega\times {\bf R}\to {\bf R}$ such that
$$\sup_{(x,\xi)\in\Omega\times {\bf R}}{{|f(x,\xi)|}\over {1+|\xi|^p}}<+\infty \eqno{(3.1)}$$
for some $p\in \left ]0, {{n+2}\over {n-2}}\right [$.\par
\smallskip
Moreover, denote by $\tilde{\cal A}$ the class of all Carath\'eodory functions $g:\Omega\times {\bf R}\to {\bf R}$ such that
$$\sup_{(x,\xi)\in\Omega\times {\bf R}}{{|g(x,\xi)|}\over {1+|\xi|^q}}<+\infty \eqno{(3.2)}$$
for some $q\in \left ] 0, {{2}\over {n-2}}\right [$.
\smallskip
Furthermore, denote by $\hat {\cal A}$ the class of all functions $h:\Omega\times {\bf R}\to {\bf R}$ of the type
$$h(x,\xi)=f(x,\xi)+\alpha(x)g(x,\xi)$$
with $f\in {\cal A}, g\in\tilde{\cal A}$ and $\alpha\in L^2(\Omega)$.
For each $h\in\hat{\cal A}$, define the functional $I_h:H^1_0(\Omega)\to {\bf R}$, by putting
$$I_h(u)=\int_{\Omega}H(x,u(x))dx$$
for all $u\in H^1_0(\Omega)$, where
$$H(x,\xi)=\int_0^{\xi}h(x,t)dt$$
for all $(x,\xi)\in \Omega\times {\bf R}$.\par
\smallskip
By classical results (involving the Sobolev embedding theorem), the functional $I_h$
turns out to be sequentially weakly continuous, of class $C^1$, with compact derivative given by
$$I_h'(u)(w)=\int_{\Omega}h(x,u(x))w(x)dx$$
for all $u,w\in H^1_0(\Omega)$.\par
\smallskip
Now, recall that, given $h\in\hat{\cal A}$, a weak solution of the problem
$$\cases {-\left ( a+b\int_{\Omega}|\nabla u(x)|^2dx\right )\Delta u =h(x,u) & in $\Omega$\cr
& \cr u=0 & on $\partial\Omega$\cr}$$
is any $u\in H^1_0(\Omega)$ such that
$$\left ( a+b\int_{\Omega}|\nabla u(x)|^2dx\right )\int_{\Omega}\nabla u(x)\nabla w(x)dx=
\int_{\Omega}h(x,u(x))w(x)$$
for all $w\in H^1_0(\Omega)$. Let $\Phi:H^1_0(\Omega)\to {\bf R}$ be the functional defined by
$$\Phi(u)={{a}\over {2}}\|u\|^2+{{b}\over {4}}\|u\|^4$$
for all $u\in H^1_0(\Omega)$.\par
\smallskip
Hence, the weak solutions of the problem are precisely the critical points in
$H^1_0(\Omega)$ of the functional $\Phi-I_h$ which is said to be the energy functional of the problem.\par
\smallskip
The central result is as follows:\par
\medskip
THEOREM 3.1. - {\it Let $n\geq 4$, let $f\in {\cal A}$ and let $g\in\tilde{\cal A}$ be such
that the set
$$\left \{x\in \Omega : \sup_{\xi\in {\bf R}}|g(x,\xi)|>0\right\}$$
has a positive measure.\par
Then,  there exist $\lambda^*\geq 0$ such that, for each $\lambda>\lambda^*$ and each 
convex set $C\subseteq L^2(\Omega)$ whose closure in $L^2(\Omega)$
contains the set $\{G(\cdot,u(\cdot)) : u\in H^1_0(\Omega)\}$,
there exists $v^*\in C$ such
that the problem
$$\cases {-\left ( a+b\int_{\Omega}|\nabla u(x)|^2dx\right )\Delta u =f(x,u)+\lambda(G(x,u)-v^*(x))g(x,u) & in $\Omega$\cr
& \cr u=0 & on $\partial\Omega$\cr}$$
has at least three weak solutions, two of which are global minima in $H^1_0(\Omega)$ of the functional
$$u\to {{a}\over {2}}\int_{\Omega}|\nabla u(x)|^2dx+{{b}\over {4}}\left ( \int_{\Omega}|\nabla u(x)|^2dx\right ) ^2-\int_{\Omega}F(x,u(x))dx
-{{\lambda}\over {2}}\int_{\Omega}|G(x,u(x))-v^*(x)|^2dx\ .$$
If, in addition,  the functional
$$u\to {{a}\over {2}}\int_{\Omega}|\nabla u(x)|^2dx+{{b}\over {4}}\left ( \int_{\Omega}|\nabla u(x)|^2dx\right ) ^2-\int_{\Omega}F(x,u(x))dx$$
has at least two global minima in $H^1_0(\Omega)$ and the function $G(x,\cdot)$ is strictly monotone for all $x\in \Omega$,
then $\lambda^*=0$.}\par
\medskip
The main tool we use to prove Theorem 3.1 is Theorem 3.C below which, in turn, is a direct consequence of
two other results recently established in [24].\par
\medskip
To state Theorem 3.C in a compact form, we now introduce some notations.\par
\smallskip
Here and in what follows, $X$ is a non-empty set, $V, Y$ are two topological
spaces, $y_0$ is a point in $Y$.\par
\smallskip
We denote by ${\cal G}$
the family of all lower semicontinuous functions $\varphi:Y\to [0,+\infty[$,
with $\varphi^{-1}(0)=\{y_0\}$, such that, for each neighbourhood $U$ of $y_0$,
one has
$$\inf_{Y\setminus U}\varphi>0\ .$$
Moreover, denote by ${\cal H}$ the family of all functions
$\Psi:X\times V\to Y$ such that, for each $x\in X$, 
$\Psi(x,\cdot)$ is continuous, injective, open, 
takes the value $y_0$ at a point $v_x$ and the function
$x\to v_x$ is not constant. 
Furthermore, denote by ${\cal M}$ the family of all
functions $J:X\to {\bf R}$ whose set of all global minima
(noted by $M_{J}$) is non-empty.\par
\smallskip
Finally, for each $\varphi\in {\cal G}$, $\Psi\in {\cal H}$ 
 and $J\in {\cal M}$, put
$$\theta(\varphi,\Psi,J)=\inf\left \{
{{J(x)-J(u)}\over {\varphi(\Psi(x,v_u))}} : (u,x)\in M_{J}\times
X\hskip 3pt \hbox {\rm with}\hskip 3pt v_x\neq v_u\right
\}\ .$$
When $X$ is a topological space, a function $\psi:X\to {\bf R}$ is said to be 
inf-compact if $\psi^{-1}(]-\infty,r])$ is compact for all $r\in {\bf R}$.
\medskip
THEOREM 3.A ([24], Theorem 3.1). - {\it Let $\varphi\in {\cal G}$, $\Psi\in {\cal H}$  
 and $J\in {\cal M}$. \par
Then, for each $\lambda>\theta(\varphi,\Psi,J)$, one has
$$\sup_{v\in V}\inf_{x\in X}
(J(x)-\lambda\varphi(\Psi(x,v)))<
\inf_{x\in X}\sup_{z\in X}
(J(x)-\lambda\varphi(\Psi(x,v_z)))\ .$$}\par
\medskip
THEOREM 3.B ([24], Theorem 3.2). - {\it Let $X$ be a topological space, $E$ a real Hausdorff
topological vector space, $C\subseteq
E$ a convex set, $f:X\times C\to
{\bf R}$ a function which is lower semicontinuous and inf-compact in $X$,
and upper semicontinuous and concave in $C$. Assume also that
$$\sup_{v\in C}\inf_{x\in X}f(x,v)<
\inf_{x\in X}\sup_{v\in C}f(x,v)\ .$$
Then, there exists $v^*\in C$ such that
the function $f(\cdot,v^*)$ has at least two global minima.}\par
\medskip
THEOREM 3.C. - {\it Let $\varphi\in {\cal G}$, $\Psi\in {\cal H}$  
 and $J\in {\cal M}$. 
Moreover, assume that $X$ is a topological space, that
$V$ is a real Hausdorff topological vector space and that
$\varphi(\Psi(x,\cdot))$ is convex and continuous for each $x\in X$. Finally, let
$\lambda>\theta(\varphi,\Psi,J)$  and let $C\subseteq V$ a convex set, with
$\{v_x : x\in X\}\subseteq \overline {C}$,
such that
 the function $x\to J(x)-\lambda\varphi(\Psi(x,v))$ is
lower semicontinuous and inf-compact in $X$ for all $v\in C$.\par
Under such hypotheses, 
there exist $v^*\in C$
 such
that the function $x\to J(x)-\lambda\varphi(\Psi(x,v^*))$ 
has at least two global minima in $X$.}\par
\smallskip
PROOF. Set
$$D=\{v_x : x\in X\}$$
and, for each $(x,v)\in X\times V$, put
$$f(x,v)=J(x)-\lambda\varphi(\Psi(x,v))\ .$$
Theorem 3.A ensures that
$$\sup_{v\in V}\inf_{x\in X}f(x,v)<
\inf_{x\in X}\sup_{v\in D}f(x,v)\ .\eqno{(3.3)}$$
But, since $f(x,\cdot)$ is continuous and $D\subseteq \overline {C}$,
 we have
$$\sup_{v\in D}f(x,v)=\sup_{v\in \overline {D}}f(x,v)\leq \sup_{v\in \overline {C}}f(x,v)
= \sup_{v\in C}f(x,v)$$
for all $x\in X$, and hence, from $(3.3)$, it follows that 
$$\sup_{v\in C}\inf_{x\in X}f(x,v)<
\inf_{x\in X}\sup_{v\in D}f(x,v)\leq \inf_{x\in X}\sup_{v\in C}f(x,v)\ .$$
At this point, the conclusion follows applying Theorem 3.B to the restriction of
the function $f$ to $X\times C$.\hfill $\bigtriangleup$\par
\medskip
{\it Proof of Theorem 3.1.} For each $\lambda\geq 0$,
$v\in L^2(\Omega)$, consider the function $h_{\lambda,v}:\Omega\times {\bf R}\to {\bf R}$
defined by 
$$h_{\lambda,v}(x,\xi)=f(x,\xi)+\lambda(G(x,\xi)-v(x))g(x,\xi)$$ 
for all $(x,\xi)\in\Omega\times {\bf R}$.
Clearly, the function $h_{\lambda,v}$ lies in
$\hat{\cal A}$ and
$$H_{\lambda,v}(x,\xi)=F(x,\xi)+{{\lambda}\over {2}}\left ( |G(x,\xi)-v(x)|^2-|v(x)|^2\right )\ .$$
 So, the weak solutions of the problem are precisely the critical points in $H^1_0(\Omega)$
of the functional $\Phi-I_{h_{\lambda,v}}$. Moreover, if $p\in \left] 0,{{n+2}\over {n-2}}\right [$ and $q\in
\left ] 0,{{2}\over {n-2}}\right [$ are such that $(3.1)$ and $(3.2)$ hold, 
 for some constant $c_{\lambda,v}$, we have
$$\int_{\Omega}|H_{\lambda,v}(x,u(x))|dx\leq c_{\lambda,v}\left ( \int_{\Omega}|u(x)|^{p+1}+
\int_{\Omega}|u(x)|^{2(q+1)}dx+1\right )$$
for all $u\in H^1_0(\Omega)$. Therefore, by the Sobolev embedding theorem, for a constant $\tilde c_{\lambda,v}$,
we have
$$\Phi(u)-I_{h_{\lambda,v}}(u)\geq {{b}\over {4}}\|u\|^4-\tilde c_{\lambda,v}(\|u\|^{p+1}+\|u\|^{2(q+1)}+1) \eqno{(3.4)}$$
for all $u\in H^1_0(\Omega)$. On the other hand, since $n\geq 4$, one has
$$\max\{p+1, 2(q+1)\}<{{2n}\over {n-2}}\leq 4\ .$$
Consequently, from $(3.4)$, we infer that
$$\lim_{\|u\|\to +\infty}(\Phi(u)-I_{h_{\lambda,v}}(u))=+\infty\ .\eqno{(3.5)}$$
Since the functional $\Phi-I_{h_{\lambda,v}}$ is sequentially weakly lower semicontinuous, by the Eberlein-Smulyan theorem and
by $(3.5)$, it follows that it is inf-weakly compact.\par
Now, we are going to apply Theorem 3.C taking $X=H^1_0(\Omega)$ with
the weak topology and $V=Y=L^2(\Omega)$ with the strong topology, and $y_0=0$. Also, we take
$$\varphi(w)={{1}\over {2}}\int_{\Omega}|w(x)|^2dx$$
for all $w\in L^2(\Omega)$. Clearly, $\varphi\in {\cal G}$. Furthermore, we take
$$\Psi(u,v)(x)=G(x,u(x))-v(x)$$
for all $u\in H^1_0(\Omega), v\in L^2(\Omega), x\in \Omega$. Clearly, $\Psi(u,v)\in L^2(\Omega)$,
 $\Psi(u,\cdot)$ is a homeomorphism, and we have
$$v_u(x)=G(x,u(x))\ .$$
We show that the map $u\to v_u$ is not constant in $H^1_0(\Omega)$. For each $x\in \Omega$,
set
$$\alpha(x)=\inf_{\xi\in {\bf R}} G(x,\xi)$$
and
$$\beta(x)=\sup_{\xi\in {\bf R}} G(x,\xi)\ .$$
Since $G$ is a Carath\'eodory is continuous, we have
$$\alpha(x)=\inf_{\xi\in {\bf Q}} G(x,\xi)$$
and
$$\beta(x)=\sup_{\xi\in {\bf Q}} G(x,\xi)\ ,$$
and so the functions $\alpha, \beta$ are measurable. 
Set
$$A=\{x\in \Omega : \alpha(x)<\beta(x)\}\ .$$
Clearly, we have
$$A=\left \{x\in \Omega : \sup_{\xi\in {\bf R}}|g(x,\xi)|>0\right\}\ .$$
Hence, by assumption, meas$(A)>0$. Then, by the classical Scorza-Dragoni theorem ([6], Theorem 2.5.19),
there exists a compact set $K\subset A$, of positive measure, such that the restriction
of $G$ to $K\times {\bf R}$ is continuous. Fix a point $\tilde x\in K$ such that
the intersection of $K$ and any ball centered at $\tilde x$ has a positive measure.
Next, fix $\xi_1, \xi_2\in {\bf R}$
such that
$$G(\tilde x,\xi_1)<G(\tilde x,\xi_2)\ .$$
By continuity, there is a closed ball $B(\tilde x,r)$ such that
$$G(x,\xi_1)<G(x,\xi_2)$$
for all $x\in K\cap B(\tilde x,r)$. Finally, consider two functions $u_1, u_2\in H^1_0(\Omega)$
which are constant in $K\cap B(\tilde x,r)$. So, we have
$$G(x,u_1(x))<G(x,u_2(x))$$
for all $x\in K\cap B(\tilde x,r)$. Hence, as meas($K\cap B(\tilde x,r))>0$, we infer that
$v_{u_1}\neq v_{u_2}$, as claimed. As a consequence, $\Psi\in {\cal H}$.  Of course,
$\varphi(\Psi(u,\cdot))$  is continuous and convex for all $u\in X$. 
Finally, take
$$J=\Phi-I_f\ .$$
Clearly, $J\in {\cal M}$. So,  for what seen above, all the assumptions of Theorem 3.C are satisfied.
Consequently, if we take 
$$\lambda^*=\theta(\varphi, \Psi,J)\eqno{(3.6)}$$ and
fix $\lambda>\lambda^*$ and a convex set  $C\subseteq L^2(\Omega)$ whose closure in $L^2(\Omega)$
contains the set $\{G(\cdot,u(\cdot)) : u\in H^1_0(\Omega)\}$, there exists $v^*\in C$ such that
the functional $\Phi-I_{h_{\lambda,v^*}}$ has at least two global minima in $H^1_0(\Omega)$ which are, therefore,
weak solutions of the problem. To guarantee the existence of a third solution, 
  denote by $k$ the inverse of the restriction of the function $at+bt^3$ to $[0,+\infty[$. 
Let $T:X\to X$ be the operator defined by
$$T(w)=\cases {{{k(\|w\|)}\over {\|w\|}}w & if $w\neq 0$\cr & \cr
0 & if $w=0$\ ,\cr}$$
 Since $k$ is continuous
and $k(0)=0$, the operator $T$ is continuous in $X$. For each $u\in X\setminus \{0\}$,
we have
$$T(\Phi'(u))=T((a+b\|u\|^2)u)={{k((a+b\|u\|^2)\|u\|)}\over {(a+b\|u\|^2)\|u\|}}(a+b\|u\|^2)u={{\|u\|}\over {(a+b\|u\|^2)\|u\|}}(a+b\|u\|^2)u=u\ .$$
In other words, $T$ is a continuous inverse of $\Phi'$. Then, since $I_{h_{\lambda,v^*}}'$ is compact, the functional
$\Phi-I_{h_{\lambda,v^*}}$ satisfies the Palais-Smale condition ([30], Example 38.25) and hence the existence of a third critical
point of the same functional is assured by Corollary 1 of [20].\par
Finally, assume that the functional $\Phi-I_f$ has at least two global minima, say $\hat u_1, \hat u_2$. Then,
the set $D:=\{x\in \Omega : \hat u_1(x)\neq \hat u_2(x)\}$ has a positive measure. By assumption, we have
$$G(x,\hat u_1(x))\neq G(x,u_2(x))$$
for all $x\in D$, and so $v_{\hat u_1}\neq v_{\hat u_2}$. Then, by definition, we have
$$0\leq\theta(\varphi,\Psi,J)\leq {{J(\hat u_1)-J(\hat u_2)}\over {\varphi(\Psi(\hat u_1,v_{\hat u_2}))}}=0$$
and so $\lambda^*=0$ in view of $(3.6)$.\hfill $\bigtriangleup$\par
\medskip
Notice the following corollary of Theorem 3.1:\par
\medskip
COROLLARY 3.1. - {\it Let $n\geq 4$, let $\nu\in {\bf R}$ and let $p\in \left ] 0,{{n+2}\over {n-2}}\right [$.\par
Then, for each $\lambda>0$ large enough and for each convex set $C\subseteq L^2(\Omega)$ whose closure in $L^2(\Omega)$
contains $H^1_0(\Omega)$, there exists $v^*\in C$ such that the problem
$$\cases {-\left ( a+b\int_{\Omega}|\nabla u(x)|^2dx\right )\Delta u =\nu|u|^{p-1}u+\lambda(u-v^*(x)) & in $\Omega$\cr
& \cr u=0 & on $\partial\Omega$\cr}$$
has at least three solutions, two of which are global minima in $H^1_0(\Omega)$ of the functional
$$u\to {{a}\over {2}}\int_{\Omega}|\nabla u(x)|^2dx+{{b}\over {4}}\left ( \int_{\Omega}|\nabla u(x)|^2dx\right ) ^2-{{\nu}\over
{p+1}}\int_{\Omega}|u(x)|^{p+1}dx
-{{\lambda}\over {2}}\int_{\Omega}|u(x)-v^*(x)|^2dx\ .$$}
\smallskip
PROOF. Apply Theorem 3.1 taking
$f(x,\xi)=|\xi|^{p-1}\xi$ and $g(x,\xi)=1$.\hfill $\bigtriangleup$\par
\medskip
REMARK 3.1. - In Theorem 3.1,  the assumption made on $g$ (besides $g\in \tilde{\cal A}$) is essential. Indeed, if $g=0$, for
$f=0$ (which is an allowed choice), the problem would have the zero solution only.\par
\medskip
REMARK 3.2. - The assumption $n\geq 4$ is likewise essential. Indeed, Corollary 3.1 does not hold if $n=3$. To see this, 
 take $p=4$ (which, when $n=3$, is compatible
with the condition $p<{{n+2}\over {n-2}}$) and observe
that the corresponding energy functional is unbounded below.\par
\medskip
Besides Corollary 3.1, among the consequences of Theorem 3.1, we highlight the following\par
\medskip
THEOREM 3.2. - {\it Let $n\geq 4$, let $f\in {\cal A}$ and let $g\in\tilde{\cal A}$ be such
the set
$$\left\{x\in \Omega : \sup_{\xi\in {\bf R}}F(x,\xi)>0\right\}$$
has a positive measure.
 Moreover, assume that, for each $x\in \Omega$, $f(x,\cdot)$ is
odd, $g(x,\cdot)$ is even and $G(x,\cdot)$ is strictly monotone.\par
Then,  for each $\lambda>0$, there exists $\mu^*> 0$ such that, for each $\mu>\mu^*$ and for each convex set $C\subseteq L^2(\Omega)$ whose closure
in $L^2(\Omega)$ contains the set $\{G(\cdot,u(\cdot)) : u\in H^1_0(\Omega)\}$, 
there exists $v^*\in C$ such
that the problem
$$\cases {-\left ( a+b\int_{\Omega}|\nabla u(x)|^2dx\right )\Delta u =\mu f(x,u)-\lambda v^*(x)g(x,u) & in $\Omega$\cr
& \cr u=0 & on $\partial\Omega$\cr}$$
has at least three weak solutions, two of which are global minima in $H^1_0(\Omega)$ of the functional
$$u\to {{a}\over {2}}\int_{\Omega}|\nabla u(x)|^2dx+{{b}\over {4}}\left ( \int_{\Omega}|\nabla u(x)|^2dx\right ) ^2-\mu\int_{\Omega}F(x,u(x))dx
+\lambda\int_{\Omega}v^*(x)G(x,u(x))dx\ .$$}\par
\smallskip
PROOF. Set
$$D=\left \{x\in \Omega : \sup_{\xi\in {\bf R}}F(x,\xi)>0\right\}\ .$$
By assumption, meas$(D)>0$. Then, by  the Scorza-Dragoni theorem,
there exists a compact set $K\subset D$, of positive measure, such that the restriction
of $F$ to $K\times {\bf R}$ is continuous. Fix a point $\hat x\in K$ such that
the intersection of $K$ and any ball centered at $\hat x$ has a positive measure.
Choose $\hat \xi\in {\bf R}$ so that $F(\hat x,\hat\xi)>0$. By continuity, there
is $r>0$ such that
$$F(x,\hat\xi)>0$$
for all $x\in K\cap B(\hat x,r)$. Set
$$M=\sup_{(x,\xi)\in \Omega\times [-|\hat \xi|,|\hat \xi|]}|F(x,\xi)|\ .$$
Since $f\in {\cal A}$, we have $M<+\infty$. Next, choose an open set $\tilde\Omega$
such that
$$K\cap B(\hat x,r)\subset\tilde\Omega\subset\Omega$$
and
$$\hbox {\rm meas}(\tilde\Omega\setminus (K\cap B(\hat x,r))<{{\int_{K\cap B(\hat x,r)}F(x,\hat \xi)dx}\over {M}}\ .$$
Finally, choose a function $\tilde u\in H^1_0(\Omega)$ such that
$$\tilde u(x)=\hat\xi$$
for all $x\in K\cap B(x,r)$,
$$\tilde u(x)=0$$
for all $x\in \Omega\setminus\tilde\Omega$ and
$$|\tilde u(x)|\leq |\hat\xi|$$
for all $x\in\Omega$. Thus, we have
$$\int_{\Omega}F(x,\tilde u(x))dx=\int_{K\cap B(\hat x,r)}F(x,\hat \xi)dx+\int_{\tilde\Omega\setminus (K\cap B(\hat x,r)}F(x,\tilde u(x))dx$$
$$>\int_{K\cap B(\hat x,r)}F(x,\hat \xi)dx-M\hbox {\rm meas}(\tilde\Omega\setminus (K\cap B(\hat x,r))>0\ .$$
Now, fix any $\lambda>0$ and set
$$\mu^*={{\Phi(\tilde u)+{{\lambda}\over {2}}I_{Gg}(\tilde u)}\over {I_f(\tilde u)}}\ .$$
Fix $\mu>\mu^*$. Hence
$$\Phi(\tilde u)-\mu I_f(\tilde u)+{{\lambda}\over {2}}I_{Gg}(\tilde u)<0\ .$$
From this, we infer that the functional $\Phi-\mu I_f+{{\lambda}\over {2}}I_{Gg}$ possesses at least to global
minima since it is even. At this point, we can apply Theorem 3.1 to the functions $g$ and $\mu f-\lambda Gg$. Our current
conclusion follows from the one of Theorem 3.1
since we have $\lambda^*=0$ and hence we can take the same fixed $\lambda>0$.\hfill $\bigtriangleup$
\bigskip
{\bf 4. A conjecture on disconnectedness versus infinitely many solutions}\par
\bigskip
CONJECTURE 4.1. - {\it Let $\Omega\subset {\bf R}^n$ be a smooth bounded domain, with $n\geq 3$. Let $\tau$ be the
strongest vector topology on $H^1_0(\Omega)$.\par
Then, there exists a continuous function $f:{\bf R}\to {\bf R}$, with
$$\sup_{\xi\in {\bf R}}{{|f(\xi)|}\over {1+|\xi|^{{n+2}\over {n-2}}}}<+\infty\ ,$$
such that the set
$$\left \{ (u,v)\in H^1_0(\Omega)\times  H^1_0(\Omega) : \int_{\Omega}\nabla u(x)\nabla v(x)dx-\int_{\Omega}f(u(x))v(x)dx=1
\right\}$$
is disconnected in $(H^1_0(\Omega),\tau)\times ( H^1_0(\Omega),\tau)$.}\par
\medskip
The importance of Conjecture 4.1 is shown by Proposition 4.3 below. But, first the relevant theory should be fixed.\par
\smallskip
The central abstract result, obtained in [21], is as follows (see also [16]):\par
\medskip
THEOREM 4.1. - {\it Let $X$ be a connected topological space, let
$E$ be a real topological vector space, with
 topological dual $E^{*}$, and let $A:X\rightarrow E^{*}$ be an operator
such that the set
$$\{y\in E : x\rightarrow \langle A(x),y\rangle
\hskip 5pt is\hskip 5pt continuous\}$$
is dense in $E$ and the set
$$\{(x,y)\in X\times E : \langle A(x),y\rangle=1\}$$
is disconnected.\par
Then, $A$ does vanish at some point of $X$.}\par
\smallskip
PROOF. Denote by $p_{X}$ the projection from $X\times E$ onto
$X$. Moreover, for any $C\subseteq X\times E$, $x\in X$, put
$$C_{x}=\{y\in E : (x,y)\in C\}.$$
Arguing by contradiction, assume that $A(x)\neq 0$ for all $x\in X$.
Denote by $\Gamma$ the set
$$\{(x,y)\in X\times E :
\langle A(x),y\rangle=1\}.$$
Since $\Gamma$
is disconnected,
 there are two open sets $\Omega_{1}, \Omega_{2}\subseteq X\times 
E$
such that $$\Omega_{1}\cap \Gamma\neq \emptyset,\hskip 3pt \Omega_{2}
\cap \Gamma\neq \emptyset,\hskip 3pt \Omega_{1}\cap\Omega_{2}\cap\Gamma=
\emptyset,\hskip 3pt \Gamma\subseteq \Omega_{1}\cup \Omega_{2}.$$
  We now prove that
  $p_{X}(\Omega_{1}\cap \Gamma)$  is open in $X$. So, let
$(x_{0},y_{0})\in \Omega_{1}\cap \Gamma$. Since $E$ is locally
connected ([28], p.35), there are 
  a neighbourhood $U_{0}$
of $x_{0}$ in $X$ and an open connected 
  neighbourhood $V_{0}$ of $y_{0}$ in $E$ such
that $U_{0}\times V_{0}\subseteq \Omega_{1}$. 
 Since 
$\langle A(x_{0}),\cdot\rangle$ is a non-null
 continuous linear functional, it has no
local extrema. Consequently, since $\langle A(x_{0}),y_{0}\rangle=1$,
 the sets
 $$\{u\in V_{0} : \langle A(x_{0}),u\rangle<1\},$$
$$\{u\in V_{0} : \langle A(x_{0}),u\rangle>1\}$$ are both non-empty and
open. Then, thanks to our density assumption, there are $u_{1}, u_{2}\in
V_{0}$ such that the
set $$\{x\in U_{0} : \langle A(x),u_{1}\rangle<1
 <\langle A(x),u_{2}\rangle\}$$ is a
neighbourhood of $x_{0}$. Then, if $x$ is in this set, due to the
connectedness of $V_{0}$, there is some $y\in V_{0}$ such that
$\langle A(x),y\rangle=1$, and so, $x$ actually lies in
 $p_{X}(\Omega_{1}\cap \Gamma)$, as desired. Likewise, it is seen
that $p_{X}(\Omega_{2}\cap \Gamma)$ is open. Now,
observe that, for any $x\in X$, the set $\{x\}\times \Gamma_{x}$ is
non-empty and connected, and so it is
contained either in $\Omega_{1}$ or in $\Omega_{2}$.
 Summarizing, we then have that the sets
$p_{X}(\Omega_{1}\cap \Gamma)$ and $p_{X}(\Omega_{2}\cap \Gamma)$
are non-empty, open, disjoint and cover $X$. 
 Hence,
$X$ would be disconnected, a contradiction.\hfill $\bigtriangleup$
\par
\smallskip
Once Theorem 4.1 has been obtained, we can state the following formally
more complete result:\par
\medskip
THEOREM 4.2. - {\it Let $X$ be a topological space, let
$E$ be a real topological vector space, and let
 $A:X\rightarrow E^{*}$
be such that the set
$$\{y\in E : x\rightarrow \langle A(x),y\rangle
\hskip 5pt is\hskip 5pt continuous\}$$
is dense in $E$.\par
Then, the following assertions are equivalent}:\par
\noindent
(i)\hskip 5pt {\it The set
 $$\{(x,y)\in X\times E : \langle A(x),y\rangle=1\}$$
is disconnected.}\par
\noindent
(ii)\hskip 5pt {\it The set $X\setminus A^{-1}(0)$ is disconnected.}\par
\smallskip
PROOF. Let (i) hold. Since $$\{(x,y)\in X\times E : \langle
 A(x),y\rangle=1\}=
\{(x,y)\in (X\setminus A^{-1}(0))\times E : \langle
 A(x),y\rangle=1\},$$ if  
$X\setminus A^{-1}(0)$ were connected, we could apply Theorem 4.1
to $A_{|(X\setminus A^{-1}(0))}$, and so $A$ would vanish at some
point of $X\setminus A^{-1}(0)$, which is absurd.\par
Conversely, if (ii) holds, then (i) follows at once observing
that, with the notations of the proof of Theorem 4.1, one has
 $X\setminus A^{-1}(0)=p_{X}(\Gamma).$\hfill $\bigtriangleup$\par
\medskip
REMARK 4.1. - When $X$ is a connected topological space, $E$ is an
 infinite-dimensional
real vector space (with algebraic dual $E'$), and $A:X\rightarrow E'$ is
a $\sigma(E',E)$-continuous operator, one could try to apply
Theorem 4.1 endowing $E$ with the strongest vector topology 
([15], p.53).\par
\medskip
REMARK 4.2. - In Theorem 4.1, the role of the constant $1$ can actually
be assumed by any continuous real function on $X$. Precisely, we have
the following\par
\medskip
PROPOSITION 4.1. - {\it Let $X$ be a topological space, let $E$ be
a real topological vector space, and let $A:X\rightarrow E'$.
Assume that, for some continuous function $\alpha : X\rightarrow
{\bf R}$, the set
$$\Lambda:=\{(x,y)\in X\times E :\langle A(x),y\rangle=\alpha(x)\}$$
is disconnected.\par
Then, either $A(x)=0$ for some $x\in X$, or the set
$$\Gamma:=\{(x,y)\in X\times E :\langle A(x),y\rangle=1\}$$
is disconnected.}\par
\smallskip
PROOF. Assume that $A^{-1}(0)=\emptyset$. So, $p_{X}(\Gamma)=X$.
 Consider the function $f:X\times E\rightarrow X\times E$
defined by putting
 $f(x,y)=(x,\alpha(x)y)$ for all $(x,y)\in X\times E$. Of course,
$f$ is continuous.
 Arguing by contradiction, assume that $\Gamma$ is connected. Then,
$f(\Gamma)$ is connected too. Now, observe that
$$\Lambda=\bigcup_{x\in \alpha^{-1}(0)}(f(\Gamma)\cup (\{x\}\times
\Lambda_{x})).$$
Furthermore, note that, if $x\in \alpha^{-1}(0)$, then
$(x,0)\in f(\Gamma)\cap (\{x\}\times \Lambda_{x})$,
and so $f(\Gamma)\cup (\{x\}\times \Lambda_{x})$ is connected. In turn,
the sets $f(\Gamma)\cup (\{x\}\times \Lambda_{x})$ ($x\in \alpha^{-1}(0))$
are clearly pairwise non-disjoint, and hence
 $\Lambda$ is connected, a contradiction.\hfill $\bigtriangleup$\par
\medskip
In [22], the following proposition was pointed out:\par
\medskip
PROPOSITION 4.2 ([22], Proposition 3). - {\it Let $E$ be an infinite-dimensional Hausdorff
topological vector space and $K$ 
a relatively compact subset of $E$.\par
Then, the set $E\setminus K$ is connected.}\par
\medskip
Finally, as said, the following proposition shows the importance of Conjecture 4.1:\par
\medskip
PROPOSITION 4.3. - {\it Let $f$ be a function satisfying Conjecture 4.1.\par
Then, the problem
$$\cases {-\Delta u=f(u) & in $\Omega$\cr & \cr
u=0 & on $\partial \Omega$\cr}$$
has infinitely many weak solutions.}\par
\smallskip
PROOF. 
Let $X=W^{1,2}_{0}(\Omega)$, with the usual norm $\|u\|=(\int_{\Omega}
|\nabla u(x)|^2dx)^{1\over 2}$.
For $0<q\leq {{n+2}\over {n-2}}$ and $f\in {\cal A}_{q}$, put
$$J(u)={{1}\over {2}}\int_{\Omega}|\nabla u(x)|^{2}dx-
\int_{\Omega}\left ( \int_{0}^{u(x)}f(\xi)d\xi\right ) dx$$
for all $u\in X$.\par
\smallskip
So, the functional $J$ is of class $C^1$ on
$X$ and one has
$$J'(u)(v)=\int_{\Omega}\nabla u(x)\nabla v(x)dx-
\int_{\Omega}f(x,u(x))v(x)dx$$
for all $u, v\in X$. Hence, the critical points of $J$ in $X$ are exactly
the weak solutions of the problem.
Since $J$ is of class $C^1$, clearly
the operator $J':X\to X^{*}$ is $\tau$-weakly-star continuous. Hence, by
Theorem 4.2, the set $X\setminus
(J')^{-1}(0)$ is $\tau$-disconnected. Then, due to Proposition 4.2, 
 the set $(J')^{-1}(0)$ is not
$\tau$-relatively compact, and hence is infinite.\hfill $\bigtriangleup$\par

\vfill\eject
\centerline {\bf References}\par
\bigskip
\bigskip
\noindent
[1]\hskip 5pt A. BAHRI and J.-M. CORON, {\it Sur une \'equation elliptique non lin\'eaire avec l'exposant critique de Sobolev},
C. R. Acad. Sci. Paris S\'er. I Math., {\bf 301} (1985), 345-348.\par
\smallskip
\noindent
[2]\hskip 5pt A. BAHRI and J.-M. CORON, {\it On a nonlinear elliptic equation involving the critical Sobolev exponent: the effect of the topology of the domain}, Comm. Pure Appl. Math., {\bf 41} (1988),
253-294.\par
\smallskip
\noindent
[3]\hskip 5pt M. CLAPP and F. PACELLA, {\it Multiple solutions to the pure critical exponent in domains with a hole of arbitrary syze},
Math. Z., {\bf 259} (2008), 575-589.\par
\smallskip
\noindent
[4]\hskip 5pt J.-M. CORON, {\it Topologie et cas limite des injections
de Sobolev}, C. R. Acad. Sci. Paris S\'er. I Math., {\bf 299} (1984),
209-212.\par
\smallskip
\noindent
[5]\hskip 5pt E. N. DANCER, {\it A note on an equation with critical
exponent}, Bull. London Math. Soc., {\bf 20} (1988), 600-602.\par
\smallskip
\noindent
[6]\hskip 5pt Z. DENKOWSKI, S. MIG\'ORSKI and N. S. PAPAGEORGIOU,
{\it An Introduction to Nonlinear Analysis: Theory}, Kluwer Academic
Publishers, 2003.\par
\smallskip
\noindent
[7]\hskip 5pt L. DIENING, C. KREUZER and S. SCHWARZACHER,  {\it Convex hull property and maximum principle for finite
element minimisers of general convex functionals}, Numer. Math., {\bf 124} (2013), 685-700.\par
\smallskip
\noindent
[8]\hskip 5pt R. ENGELKING, {\it Theory of dimensions, finite and infinite}, Heldermann Verlag, 1995.\par
\smallskip
\noindent
[9]\hskip 5pt X. L. FAN, {\it A remark on Ricceri's conjecture for a
class of nonlinear eigenvalue problems}, J. Math. Anal. Appl., 
{\bf 349} (2009), 436-442.\par
\smallskip
\noindent
[10]\hskip 5pt X. L. FAN, {\it On Ricceri's conjecture for a class of nonlinear eigenvalue problems}, Appl. Math. Lett., {\bf 22} (2009),
1386-1389.\par
\smallskip
\noindent
[11]\hskip 5pt X. L. FAN and B. RICCERI, {\it On the Dirichlet problem
involving nonlinearities with non-positive primitive: a problem and a
remark}, Appl. Anal., {\bf 89} (2010), 189-192.\par
\smallskip
\noindent
[12]\hskip 5pt N. HIRANO, {\it Existence of nontrivial solutions for a semilinear elliptic problem with supercritical exponent}, Nonlinear Anal., {\bf 55} (2003), 543-556.\par
\smallskip
\noindent
[13]\hskip 5pt N. I. KATZOURAKIS, {\it Maximum principles for vectorial approximate minimizers of nonconvex functionals},
 Calc. Var. Partial Differ. Equ., {\bf 46} (2013), 505-522.\par
\smallskip
\noindent
[14]\hskip 5pt J. L. KAZDAN and F. W. WARNER,  {\it Remarks on some
quasilinear elliptic equations}, Comm. Pure Appl. Math., {\bf 28}
(1975), 567-597.\par
\smallskip
\noindent
[15]\hskip 5pt J. L. KELLEY and I. NAMIOKA, {\it Linear topological spaces}, Van Nostrand, 1963.\par
\smallskip
\noindent
[16]\hskip 5pt A. J. B. LOPES-PINTO,  {\it On a new result on the existence of zeros due to Ricceri}, 
J. Convex Anal., {\bf 5} (1998), 57-62.
\smallskip
\noindent
[17]\hskip 5pt D. PASSASEO, {\it Multiplicity of positive solutions of nonlinear elliptic equations with critical Sobolev exponent in some contractible 
domains}, Manuscripta Math., {\bf 65} (1989), 147-175.\par
\smallskip
\noindent
[18]\hskip 5pt D. PASSASEO, {\it Nontrivial solutions of elliptic equations with supercritical exponent in contractible domains}, 
Duke Math. J., {\bf 92} (1998), 429-457.\par
\smallskip
\noindent
[19]\hskip 5pt S. I. POHOZAEV, {\it Eigenfunctions of the equation
$\Delta u +\lambda f(u)=0$}, Soviet Math. Dokl., {\bf 6} (1965), 1408-1411.
\smallskip
\noindent
[20]\hskip 5pt P. PUCCI and J. SERRIN, {\it A mountain pass theorem},
J. Differential Equations, {\bf 60} (1985), 142-149.\par
\smallskip
\noindent
[21]\hskip 5pt B. RICCERI, {\it Existence of zeros via disconnectedness},
J. Convex Anal., {\bf 2} (1995), 287-290.\par
\smallskip
\noindent
[22]\hskip 5pt B. RICCERI, {\it Applications of a theorem concerning sets
with connected sections}, Topol. Methods Nonlinear Anal., {\bf 5} (1995), 237-248.\par
\smallskip
\noindent
[23]\hskip 5pt B. RICCERI, {\it A remark on a class of nonlinear eigenvalue problems}, Nonlinear Anal., {\bf 69} (2008), 2964-2968.\par
\smallskip
\noindent
[24]\hskip 5pt B. RICCERI, {\it A strict minimax inequality criterion and some of its consequences}, Positivity, {\bf 16} (2012), 455-470.\par
\smallskip
\noindent
[25]\hskip 5pt B. RICCERI, {\it Energy functionals of Kirchhoff-type problems having multiple global minima}, Nonlinear Anal., {\bf 115} (2015),  
130-136.\par
\smallskip
\noindent
[26]\hskip 5pt B. RICCERI, {\it The convex hull-like property and supported images of open sets}, Ann. Funct. Anal., {\bf 7}
(2016), 150-157.\par
\smallskip
\noindent
[27]\hskip 5pt J. SAINT RAYMOND, {\it Open differentiable mappings}, Le Matematiche, {\bf 71} (2016), 197-208.\par
\smallskip
\noindent
[28] H. H. SCHAEFER, {\it Topological vector spaces}, Springer-Verlag, 1971.
\smallskip
\noindent
[29]\hskip 5pt  E. ZEIDLER, {\it Nonlinear functional analysis and its
applications}, vol. I, Springer-Verlag, 1986.\par
\smallskip
\noindent
[30]\hskip 5pt E. ZEIDLER, {\it Nonlinear functional analysis and its
applications}, vol. III, Springer-Verlag, 1985.\par
\bigskip
\bigskip
\bigskip
\bigskip
Department of Mathematics\par
University of Catania\par
Viale A. Doria 6\par
95125 Catania\par
Italy\par
{\it e-mail address}: ricceri@dmi.unict.it

\bye